\documentclass[11pt]{amsart}
\usepackage{amsfonts,latexsym,amsthm,amssymb,amsmath,amscd,euscript,tikz, tikz-cd, afterpage}
\usepackage[alphabetic, msc-links, bibtex-style, nobysame]{amsrefs}
\usepackage{stackengine}
\usepackage{accents}
\usepackage{framed}
\usepackage{xfrac,xcolor}
\usepackage[normalem]{ulem}
\usepackage[makeroom]{cancel}
\usepackage{faktor}
\usepackage{braket}
\usepackage{pgf,tikz,pgfplots}
\pgfplotsset{compat=1.17}
\usepgfplotslibrary{fillbetween} 
\usepackage{mathrsfs}
\usetikzlibrary{arrows}
\usetikzlibrary{arrows.meta}
\definecolor{wrwrwr}{rgb}{0.3803921568627451,0.3803921568627451,0.3803921568627451}
\definecolor{rvwvcq}{rgb}{0.08235294117647059,0.396078431372549,0.7529411764705882}
\definecolor{mblue}{rgb}{0.2, 0.3, 0.8}
\definecolor{morange}{rgb}{1, 0.5, 0}
\definecolor{mgreen}{rgb}{0.1, 0.4, 0.2}
\definecolor{mred}{rgb}{0.5, 0, 0}
\definecolor{ForestGreen}{RGB}{34,139,34}
\usepackage[]{hyperref}
\hypersetup{colorlinks=true,citecolor=ForestGreen,linkcolor = blue,urlcolor =black,linkbordercolor={1 0 0}}
\usepackage{float}

\usepackage{lmodern}
\usepackage{supertabular}
\usepackage{verbatim}
\usepackage{amssymb}
\usepackage{enumerate}
\usepackage{stmaryrd}
\usepackage{bbm}
\usepackage{mathtools}
\usepackage[nopatch=footnote]{microtype}
\usepackage{setspace}
\usepackage{tikz,tikz-cd}
\usetikzlibrary{matrix,calc,positioning,arrows,arrows.meta,positioning,decorations.pathreplacing,decorations.pathmorphing,patterns,knots,}

\newcommand{\la}{\langle}
\newcommand{\rg}{\rangle}

\newtheorem{theorem}{{Theorem}}[section]
\newtheorem*{theorem*}{Theorem}
\newtheorem{lemma}[theorem]{Lemma}
\newtheorem{proposition}[theorem]{Proposition}

\newtheorem{corollary}[theorem]{Corollary}
\newtheorem*{corollary*}{Corollary}

\theoremstyle{definition}

\newtheorem{remark}{Remark}
\newtheorem{assumption}[theorem]{Assumption}

\usepackage{lmodern,url,enumerate,mathtools,microtype}
\usepackage[hmargin = 1in,vmargin=1in]{geometry}
\usepackage{graphicx}
\usepackage{subcaption}

\newcommand{\ve}{\varepsilon}

\newcommand{\mr}[1]{{\rm #1}}
\mathchardef\originall=\mathcode`l
\mathcode`l="8000
\begingroup
  \lccode`~=`l
  \lowercase{\endgroup
    \def~{%
      \ifnum\mathgroup=\m@ne
        \ell             % Use \ell for variables (mathgroup -1)
      \else
        \originall       % Use normal 'l' for operators/text
      \fi
    }%
  }

\newcommand{\bR}{\mathbb{R}}

\newcommand{\nc}{\newcommand}

\nc{\on}{\operatorname}
\nc{\p}{\partial}
\nc{\ol}{\overline}
\nc{\ul}{\underline}
\nc{\pa}{\partial}

\nc{\pb}{\partial_b}
\nc{\pc}{\partial_c}
\nc{\pd}{\partial_d}
\nc{\pe}{\partial_e}
\nc{\pf}{\partial_f}
\nc{\pg}{\partial_g}
\nc{\ph}{\partial_h}
\nc{\pari}{\partiaQ_i}
\nc{\pj}{\partial_j}
\nc{\pk}{\partial_k}
\nc{\pl}{\partial_l}
\nc{\pell}{\partial_\ell}
\nc{\parm}{\partial_m}
\nc{\pn}{\partial_n}
\nc{\po}{\partial_o}
\nc{\pp}{\partial_p}
\nc{\pq}{\partial_q}
\nc{\pr}{\partial_r}
\nc{\ps}{\partial_s}
\nc{\pt}{\partial_t}
\nc{\pu}{\partial_u}
\nc{\pv}{\partial_v}
\nc{\pw}{\partial_w}
\nc{\px}{\partial_x}
\nc{\py}{\partial_x}
\nc{\pz}{\partial_z}

\nc{\Spec}{\on{Spec}}

\nc{\sn}{\mr{sn}}
\nc{\cn}{\mr{cn}}
\nc{\dn}{\mr{dn}}
\nc{\gt}{>}
\nc{\lt}{<}

\allowdisplaybreaks

\newcommand{\@dotsep}{4.5}

\def\@tocline#1#2#3#4#5#6#7{%
  \relax
  \ifnum #1>\c@tocdepth
  \else
    \par \addpenalty\@secpenalty \addvspace{#2}%
    \begingroup
      \hyphenpenalty\@M
      \@ifempty{#4}{%
        \@tempdima\csname r@tocindent\number#1\endcsname\relax
      }{%
        \@tempdima#4\relax
      }%
      \parindent\z@
      \leftskip#3\relax
      \advance\leftskip\@tempdima\relax
      \rightskip\@pnumwidth plus4em
      \parfillskip-\@pnumwidth
      #5\leavevmode\hskip-\@tempdima #6\nobreak
      \leaders\hbox{$\m@th
        \mkern \@dotsep mu.\mkern \@dotsep mu$}\hfill
      \nobreak\hbox to\@pnumwidth{\@tocpagenum{#7}}\par
    \endgroup
  \fi
}

\renewcommand{\l@section}{\@tocline{1}{2pt}{0pt}{}{\bfseries}}
\renewcommand{\l@subsection}{\@tocline{2}{0pt}{1.5em}{3.2em}{}}
\renewcommand{\l@subsubsection}{\@tocline{3}{0pt}{4.7em}{3.8em}{}}

\makeatother

\title{A positive mass theorem for manifolds with multiple corners}

\date{\today}
\author[C. Alves Sabín]{Carmen Alves Sabín}
\author[A. Arloro]{Andrea Arloro}
\author[B. Qin]{Boqian Qin}
\author[E. Wizel]{Ella Wizel}
\address{ \vspace*{-0.2in} {} \newline Department of Mathematics, Columbia University
\newline {\href{mailto:ca3104@columbia.edu}{ca3104@columbia.edu}}
\newline {\href{mailto:aa5200@columbia.edu}{aa5200@columbia.edu}}
\newline {\href{mailto:bq2170@columbia.edu}{bq2170@columbia.edu}}
\newline {\href{mailto:ew2869@columbia.edu}{ella.wizel@columbia.edu}}}

\begin{document}

\begin{abstract}
We prove a Riemannian positive mass theorem for asymptotically flat spin manifolds with multiple corners, with no assumption on the sign of the mean curvatures of the corners or the scalar curvature of the compact regions. Our hypotheses combine a pointwise comparison of the mean curvature with quantitative bounds using the distance between the corners. The proof introduces a Riccati comparison along curves and expands on the lock principle of Frenck-Hanke-Hirsch and the positive mass theorem for creased initial data of Kazaras-Khuri-Lin.
\end{abstract}

\maketitle

\section{Introduction} 

The Riemannian positive mass theorem asserts that a complete asymptotically flat manifold of non-negative scalar curvature has non-negative ADM mass.
This quantity, denoted by $m_{\textup{ADM}}$, measures the total mass–energy of a system as seen from spatial infinity.
The positive mass theorem was first established in low dimensions by Schoen-Yau~\cite{schoen-yau}, using minimal hypersurfaces, and by Witten~\cite{witten}, using spinors.
Since then, it has led to numerous developments in scalar curvature geometry and mathematical relativity.
Recently, Brendle-Wang~\cites{brendle-pmt, brendle-spacetime-pmt} developed a new dimension descent scheme to establish the Riemannian positive mass inequality and the spacetime positive energy theorem in every dimension.
We refer to~\cites{he-shi-yu, five-guys, hklz, kww} for related recent developments.

An important direction in the study of the positive mass theorem concerns its validity for metrics of lower regularity.
In particular, Miao~\cite{Miao2002} and Shi-Tam~\cite{shi-tam} established positive mass theorems for manifolds with metric singularities along hypersurfaces, assuming that the scalar curvature is non-negative away from the singular set and that the mean curvature does not increase along each hypersurface, in the direction of the asymptotically flat end.
The latter condition can be understood as the non-negativity of the singular part of the scalar curvature, and has since been expanded by various authors~\cites{lee-lefloch, tsang, KazarasKhuriLin2025}.
The latter work, due to Kazaras, Khuri, and Lin, introduced a matching condition for Bartnik data that allows hyperbolic rotations of the normal bundles along the singular hypersurface, called a \textit{crease}.

Using the positive mass theorem for creased initial data of~\cite{KazarasKhuriLin2025}, Frenck-Hanke-Hirsch~\cite{FrenckHankeHirsch2026} established a Riemannian positive mass theorem in which mean-concave corners are permitted, provided that their mean curvature defects are compensated by sufficiently mean-convex corners.
This phenomenon, which they refer to as the \textit{lock principle} for scalar curvature, establishes the non-negativity of the ADM mass even in situations where the scalar curvature fails to be non-negative in the distributional sense.
The key insight is to construct an auxiliary spacetime initial data set whose second fundamental form transfers the mean curvature defects between successive corners.

In this paper, we develop a quantitative extension of the lock principle, where the compensation of mean curvature entails the scalar curvature and the distances between the singular hypersurfaces, in addition to the jumps in mean curvature.
Our approach is motivated by the distance estimates and shielding principles developed by Cecchini-Zeidler~\cites{CecchiniZeidler2024 , cecchinilesourdzeidler2026sheilds}, Lesourd-Unger-Yau~\cite{lesourd-unger-yau}, and Cecchini-Lesourd-Zeidler~\cite{cecchinilesourdzeidler2026sheilds}.
Related ideas were also proposed in~\cite{FrenckHankeHirsch2026}*{\S 3}.

We now describe the geometric setting and state our main result.

\begin{assumption}\label{assumption}
Let $n \geq 3$ and let $M^n$ be a connected spin manifold without boundary, decomposed along finitely many disjoint closed two-sided hypersurfaces as
\[
M=M_0\cup_{\Sigma_{1}}M_1\cup_{\Sigma_2}\cdots\cup_{\Sigma_N}M_N = M_0 \cup \bigcup_{i=1}^N \left( \Sigma_i \cup \on{int}(M_i) \right), 
\]
where the $\{ M_i \}_{i=0}^{N-1}$ are compact manifolds, $M_N$ is asymptotically flat, and all pieces are connected.
We have $\partial M_0 = \Sigma_1$, $\partial M_N = \Sigma_N$, an $\partial M_i = \Sigma_i \sqcup \Sigma_{i+1}$ for $1 \leq i < N$.
We suppose that $M^n$ is equipped with a $C^0$ piecewise Riemannian metric $g$ whose restrictions $g_i := g|_{M_i}$ are smooth up to the boundary, with the same induced metric $g_{i-1}|_{\Sigma_i} = g_i|_{\Sigma_i}$ on either side of each boundary $\Sigma_i$ and extending smoothly to the corresponding boundary hypersurfaces.

On each of the finitely many ends of $M_N$, we assume that
\[
\partial^{\ell} (g_{ij} - \delta_{ij}) = O ( |x|^{-\tau - \ell}), \qquad 0 \leq \ell \leq 2, \qquad \tau > \tfrac{n-2}{2}
\]
in asymptotically Euclidean coordinates.
Moreover, we suppose that $R_{g_N} \in L^1(M_N)$.
\end{assumption}

We orient the unit normal along $\Sigma_i$ as pointing toward the asymptotically flat end, namely from $M_{i-1}$ toward $M_i$.
We define by $H_{i,-} := H(\Sigma_i \subset M_{i-1})$ and $H_{i,+} := H(\Sigma_i \subset M_i)$ the mean curvatures computed on the two sides of the hypersurface $\Sigma_i$ with respect to this normal.

To introduce our positive mass theorem, consider $p \geq 0$ and $r \geq - \frac{n}{n-1} p^2$.
For $D \geq 0$, we denote by $\Phi_r(p,D)$ the time-$D$ value of the unique solution of the initial value problem for the Riccati ODE
\begin{equation}\label{eqn:Phi-Riccati}
    \Phi_r[ p,D] :=  y(D), \qquad \text{where} \qquad y' = \tfrac{n}{2(n-1)} y^2 + \tfrac{1}{2} r, \; y(0) = p .
\end{equation}
If the solution $y$ of this initial value problem blows up in finite time $\tilde{D} \leq D$, we set $\Phi_r[p,D] = +\infty$.

Our main result is a positive mass theorem for manifolds satisfying Assumption~\ref{assumption}.

\begin{theorem}\label{thm:better-theorem}
Let $M^n$ be a manifold satisfying Assumption \ref{assumption}.
We suppose that there exist coefficients $a_i, b_i \geq 0$ such that the following conditions are satisfied.
\begin{enumerate}[(i)]
    \item We have $b_N = 0$ and $R_{g_N} \geq 0$.
    Moreover, the initial region satisfies
    \[
    R_{g_0} + \tfrac{n}{n-1} a_1^2 \geq 0 \quad \text{on } \; M_0.
    \]
    \item For each $1 \leq i < N$, we have either $b_i = a_{i+1}$ and $R_{g_i} + \frac{n}{n-1} b_i^2 \geq 0$ on $M_i$, or $b_i \neq a_{i+1}$ and there exists a real number $r_i \in \bR$ such that
    \begin{align}
        - \tfrac{n}{n-1} ( \min \{ b_i, a_{i+1} \})^2 &< r_i \leq \min_{M_i} R_{g_i}, \label{eqn:ri-condition} \\
        \max\{ b_i, a_{i+1} \} &< \Phi_{r_i} [ \min \{ b_i, a_{i+1} \} , \textup{dist} (\Sigma_i, \Sigma_{i+1}) ] + \tfrac{1}{2} \inf_{\gamma: \Sigma_i \to \Sigma_{i+1}} \int_{\gamma} (R_{g_i} - r_i) \, ds. \label{eqn:integral-Phi}
    \end{align}
    Here, the infimum is taken over all closed curves $\gamma$ connecting $\Sigma_i$ to $\Sigma_{i+1}$.
    \item For each $1 \leq i \leq N$, there is a smooth positive function $t_i$ on $\Sigma_i$ such that
    \begin{align}\label{eqn:Ai-Bi-AiBi-conditions}
        \begin{split}
    %%    t_i (H_{i,-} + a_i) - (H_{i,+} + b_i) \geq 0, \qquad H_{i,-} - a_i - t_i (H_{i,+} - b_i) &\geq 0, \\
        [t_i (H_{i,-} + a_i) - (H_{i,+} + b_i)] \cdot [ (H_{i,-} - a_i) - t_i (H_{i,+} - b_i)] &\geq t_i^{-1} |\nabla^{\Sigma_i} t_i|^2
        \end{split}
    \end{align}
    and both factors on the left-hand side are non-negative.
\end{enumerate}
Then, it holds that $m_{\textup{ADM}}(M,g) \geq 0$ on every asymptotically end.
Moreover, if $b_i \neq a_{i+1}$ for some $1 \leq i < N$, then the mass of every end is strictly positive.
\end{theorem}

We highlight that $M^n$ is not required to have positive scalar curvature.
Moreover, the corners of $M^n$ have no requirement of mean convexity, and we impose no restriction on the order of favorable and unfavorable corners.
When $b_i = a_{i+1}$, the condition $(ii)$ simply requires $R_{g_i} + \frac{n}{n-1} b_i^2 \geq 0$.
In particular, our technique compensates for negative scalar curvature on a compact piece.

The proof of Theorem~\ref{thm:better-theorem} proceeds by constructing an initial data set with second fundamental form $k_i = c_i \frac{g_i}{n-1}$ on each region $M_i$, where the function $c_i$ has prescribed values on the boundary hypersurfaces and need only be Lipschitz.
In particular, in Proposition~\ref{prop:pmt-for-creased}, we obtain an extension of the positive mass theorem for creased initial data from~\cite{KazarasKhuriLin2025}*{Theorem 1.4}, assuming only Lipschitz regularity for the $2$-tensor $k$.
The dominant energy condition becomes a scalar differential inequality for $c_i$, and we use $(ii)$ to prove the existence of such a coefficient, through a Riccati comparison and a weighted distance argument.
Moreover, condition $(iii)$ recovers the matching condition of Kazaras-Khuri-Lin~\cite{KazarasKhuriLin2025}*{\S 1} for the data in question.
Finally, taking $k_N=0$ reduces the equality case of $m_{\textup{ADM}}(M,g) \geq 0$ to the Riemannian positive mass theorem.

\medskip \noindent \textbf{Acknowledgments.}
We are grateful to Sven Hirsch and Raphael Tsiamis for suggesting this project, for their invaluable guidance, and for many helpful discussions, particularly on obtaining a more general version of Theorem ~\ref{thm:better-theorem}.
We are very thankful to Professor Simon Brendle for his generous support.
We also thank the Columbia Math Undergraduate Summer Research program, particularly George Dragomir, as well as Jane Street Capital, for providing the ideal working conditions for this project to be initiated and completed.

\medskip \noindent \textbf{AI statement.} 
No stage of this work involved LLMs.

\section{Preliminaries}\label{sec:preliminaries}

We recall that an initial data set $(M^n,g,k)$ satisfying the \textbf{dominant energy condition} is a Riemannian manifold $(M^n,g)$ endowed with a symmetric $2$-tensor $k$ such that 
\begin{equation}\label{eqn:DEC-condition}
    \mu := \tfrac{1}{2} (R_g + ( \textup{tr}_g k)^2 - |k|^2_g) \geq |J|_g, \qquad \text{where} \quad J := \textup{div}_g ( k - ( \textup{tr}_g k) g).
\end{equation}
Given a locally Lipschitz function $c \in W^{1,\infty}_{\textup{loc}}(M)$, the $2$-tensor $k = \frac{c}{n-1} g$ satisfies
\[
\textup{tr}_g k = \tfrac{n}{n-1} c, \qquad |k|^2_g = \tfrac{n}{(n-1)^2} c^2, \qquad J = - dc, \qquad \mu = \tfrac{n}{2(n-1)} c^2 + \tfrac{1}{2} R_g.
\]
Consequently, the dominant energy condition~\eqref{eqn:DEC-condition} for the $2$-tensor $k = \frac{c}{n-1
} g$ becomes equivalent to
\begin{equation}\label{eqn:dec-equivalent}
    | \nabla c|_g \leq \tfrac{1}{2} \bigl( R_g + \tfrac{n}{n-1} c^2 \bigr)
\end{equation}
since $\textup{tr}_g k = \frac{n}{n-1} c$ and $|k|^2_g = \frac{n}{(n-1)^2} c^2$.

\begin{lemma}\label{lemma:band-width-estimate}
Let $(\Omega, g)$ be a compact connected smooth Riemannian manifold with boundary $\partial \Omega = \Sigma \sqcup \tilde{\Sigma}$.
We consider numbers $0 \leq p < q$ and $r$ with the property that
\[
- \tfrac{n}{n-1} p^2 < r \leq \min_{\Omega} R_g, \qquad q < \Phi_r[p, \textup{dist}_g(\Sigma, \tilde{\Sigma})] + \tfrac{1}{2} \inf_{\gamma: \Sigma \to \tilde{\Sigma}} \int_{\gamma} (R_g - r) \, ds.
\]
Here, the infimum is taken over all curves $\gamma$ in $\Omega$ joining $\Sigma \to \tilde{\Sigma}$.

Then, there exists a smooth function $c$ on $\Omega$ such that $c=p$ in a neighborhood of $\Sigma$, $c=q$ in a neighborhood of $\tilde{\Sigma}$, and $c$ satisfies the condition~\eqref{eqn:dec-equivalent}.
\end{lemma}
\begin{proof}
We choose $\delta>0$ sufficiently small, to be determined later, such that $\frac{n}{n-1} p^2 + (r-\delta)>0$.
Moreover, we consider the solution of the initial problem for the Riccati ODE
\[
y' = \tfrac{n}{2(n-1)} y^2 + \tfrac{1}{2} (r - \delta), \qquad y(0) = p.
\]
The function $y$ is strictly increasing in its maximal interval of existence and satisfies $y(t) \geq p$.

For $\ve \in (0,\ve_0]$ sufficiently small, we let $\{ K_{\ve} \}_{\ve \in (0,\ve_0]}$ be sufficiently thin collar neighborhoods of $\Sigma$, such that $\bigcap_{\ve \in (0,\ve_0]} K_{\ve} = \Sigma$.
Given a point $x$, we then define the functions
\[
d(x) := \textup{dist}(x, K_{\ve}), \qquad w(x) := \frac{1}{2} \inf_{\gamma: K_{\ve} \to x} \int_{\gamma} (R_g - r) \, ds,
\]
where the infimum is taken over all curves $\gamma$ joining a point on $K_{\ve}$ to $x$.
In particular, we observe that $d = w = 0$ on $K_{\ve}$.
Finally, we set
\[
c_0 (x) := \min \{ q, y(d(x)) + w(x) \}.
\]
Then, our construction ensures that $d=w=0$ on $K_{\ve}$, hence $c_0 = p$ throughout this collar.
Moreover, the function $d$ is $1$-Lipschitz, so Rademacher's theorem implies that is differentiable almost everywhere, with $|\nabla d| \leq 1$.
In addition, the definition of $w$ implies that 
\begin{equation}\label{eqn:derivative-bound-d-w}
|\nabla d| \leq 1, \qquad |\nabla w| \leq \tfrac{1}{2} (R_g - r)
\end{equation}
almost everywhere.
The latter property is standard, since for nearby points $x_1, x_2$, we can concatenate curves in the definition of $w$ to obtain
\[
2 \, |w(x_1) - w(x_2)| \leq \inf_{\gamma: x_1 \to x_2} \int_{\gamma} (R_g - r) \, ds,
\]
from which~\eqref{eqn:derivative-bound-d-w} follows.
See, for example, ~\cite{docarmo}*{Ch. 3, \S 3}.

In the region $\{ p < c_0 < q \}$, we therefore obtain the inequality
\begin{equation}\label{eqn:c0-inequality}
\begin{split}
|\nabla c_0| &\leq y'(d) \, |\nabla d| + |\nabla w| \leq \tfrac{n}{2(n-1)} y(d)^2 + \tfrac{1}{2} (r-\delta) + \tfrac{1}{2} (R_g - r) \\
&\leq \tfrac{n}{2(n-1)} c_0^2 + \tfrac{1}{2} (R_g - \delta).
\end{split}
\end{equation}
Moreover, on the set $\{ c_0 = p \} \cup \{ c_0 = q \}$, the gradient $\nabla c_0$ vanishes, while
\[
\tfrac{n}{2(n-1)} c_0^2 + \tfrac{1}{2} R_g \geq \tfrac{n}{2(n-1)} p^2 + \tfrac{1}{2} r> 0
\]
due to the properties of $r$.
Thus, $c_0$ satisfies the differential inequality~\eqref{eqn:dec-equivalent} strictly.
Moreover, the construction of $c_0$ ensures that $c_0 = q$ near $\tilde{\Sigma}$ and $c_0 = p$ near $\Sigma$.

Finally, since $\Omega$ is compact, a standard argument by local mollification and a partition of unity (see, for example,~\cite{AzagraFerreraLopezMesasRangel2007}*{Theorem 1}) produces, for each $\eta>0$, a smooth function $c_{\eta}$ such that 
\[
|c_{\eta} - c| < \eta, \qquad (1-\eta) \, |\nabla c_0|  < |\nabla c_{\eta}| < (1+\eta) \, |\nabla c_0|.
\]
Combining these properties with~\eqref{eqn:c0-inequality} shows that the function $c = c_{\eta}$ satisfies the condition~\eqref{eqn:dec-equivalent} upon taking $\delta, \ve, \eta>0$ to be sufficiently small.
After decreasing $\eta$, if necessary, we can also ensure that $c=p$ near $\Sigma$ and $c=q$ near $\tilde{\Sigma}$.
This proves our assertion.
\end{proof}

\begin{remark}\label{rmk:generality}
In Lemma~\ref{lemma:band-width-estimate}, relaxing the asssumptions on $p,q,r$ to non-strict inequalities allows us to produce a Lipschitz function $c$ satisfying~\eqref{eqn:dec-equivalent} almost everywhere, with $c|_{\Sigma}=p$ and $c|_{\Sigma}=q$.
The argument employed above shows the following more general property: consider a Lipschitz function $u$ on a compact smooth Riemannian manifold-with-boundary $(\Omega,g)$, which is smooth on collar neighborhoods of $\partial \Omega$ and satisfies $|\nabla u| \leq F(x,u) - \sigma$, where $F$ is smooth and $\sigma>0$.
Then, $u$ admits smooth uniform approximations $\{ u_{\ve} \}_{ \ve \in (0,\ve_0]}$, agreeing with $u$ on smaller collar neighborhoods of $\partial \Omega$ and satisfying $|\nabla u_{\ve}| \leq F(x,u_{\ve})$.

As a particular case of Lemma~\ref{lemma:band-width-estimate}, we obtain a function $c$ with the desired properties whenever
\[
\int_p^q \frac{2 (n-1)}{nu^2 + (n-1) \min_{\Omega} R_g} \, du < \textup{dist}_g(\Sigma, \tilde{\Sigma}).
\]
The left-hand side of this expression can be computed via explicit trigonometric, rational, and hyperbolic functions upon setting $\kappa = \frac{1}{n(n-1)} \min_{\Omega} R_g$.
Moreover, if $p=q$, then the constant function $c \equiv p$ satisfies the condition~\eqref{eqn:dec-equivalent} whenever $\min_{\Omega} R_g + \frac{n}{n-1} p^2 \geq 0$.
\end{remark}

The above results can also be extended to the case $p = \min_{\Omega} R_g = 0$.
\begin{corollary}\label{cor:localization}
Let $(\Omega, g)$ be a compact Riemannian band with $R_g \geq 0$.
Suppose that there is a separating sub-band $B \subset \Omega$, with boundary components $S, \tilde{S}$, satisfying $R_g \geq r > 0$ on $B$.
If
\[
q < \Phi_r[0, \textup{dist}_B (S, \tilde{S})],
\]
then there exists a smooth function $c : \Omega \to [0,q]$ equal to zero near the first boundary of $\Omega$, equal to $q$ near the second, and satisfying the property~\eqref{eqn:dec-equivalent}.
In particular, if
\[
\textup{dist}_B(S, \tilde{S})^2 \geq \tfrac{n-1}{n} \pi^2 r^{-1},
\]
then we can find such a function $c$ for any finite $q \geq 0$.
\end{corollary}
\begin{proof}
This result follows from the argument of Lemma~\ref{lemma:band-width-estimate}.
Indeed, the strict inequality $r>0$ shows that we can first apply that result inside $B$ to produce a function $c_B$ that satisfies~\eqref{eqn:dec-equivalent}, equals $0$ near $S$, and equals $q$ near $\tilde{S}$.
We then extend $c_B$ by zero on the first side of $B$ and by $q$ on the second.
Since $R_g \geq 0$, both constant extensions satisfy~\eqref{eqn:dec-equivalent}.
Finally, we can reverse the two boundary components to also obtain a function $c$ decreasing from $q$ to zero over $\partial \Omega$.

For the second application, we observe that the solution of the ODE $y' = \frac{n}{2(n-1)} y^2 + \frac{1}{2} r$ starting at zero blows up at time $\sqrt{\frac{n-1}{nr}} \pi$.
Thus, if $\textup{dist}_B(S,\tilde{S}) \geq \frac{n-1}{n} \pi^2 r^{-1}$, our computation shows that $q < \Phi_r[0, \textup{dist}_B(S,\tilde{S})] = + \infty$ for any finite $q \geq 0$.
This proves our assertion.
\end{proof}

\subsection{The crease condition}

Next, we show that the condition $(iii)$ of Theorem~\ref{thm:better-theorem} is equivalent to the crease condition (1.7) of~\cite{KazarasKhuriLin2025}*{Definition 1.2}.
This also extends~\cite{FrenckHankeHirsch2026}*{Lemma 2}.

Let us consider an orthonormal frame $(\nu, \tau)$ with $\la \nu, \nu \rg = 1$ and $\la \tau, \tau \rg = - 1$, where $\nu$ is the unit spatial normal and $\tau$ is the future timelike normal.
For the two sides of a hypersurface, the Bartnik data associated with the initial data set $(M,g,k)$ consist of a mean curvature vector and a normal connection, given by
\begin{equation}\label{eqn:bartnik-data}
    \vec{H} = H \nu - ( \textup{tr}_{\Sigma} k) \, \tau, \qquad \beta(Z) = k(\nu, Z), \qquad \text{where } \; Z \in T \Sigma.
\end{equation}
Let $F_{\vartheta}$ be a hyperbolic rotation, namely the normal bunle isometry with hyperbolic angle $\vartheta$,
\[
F_{\vartheta} := \begin{pmatrix}
    \cosh \vartheta & - \sinh \vartheta \\
    - \sinh \vartheta & \cosh \vartheta
\end{pmatrix}.
\]
As in~\cite{KazarasKhuriLin2025}*{\S 1}, we consider the vector $X$ and the $1$-form $\beta^{\Delta}$ given by
\[
X := F_{\vartheta} \vec{H}_- - \vec{H}_+ = X_1 \nu_+ + X_2 \tau_+, \qquad \beta^{\Delta} := \beta_+ - \beta_- - df.
\]
In this setting, the crease condition~\cite{KazarasKhuriLin2025}*{(1.2)} assumes the form
\begin{equation}\label{eqn:crease-condition}
    X_1 \geq \sqrt{X_2^2 + |\beta^{\Delta}|^2}.
\end{equation}
\begin{lemma}\label{lemma:KKL}
Let $\Sigma$ be a smooth two-sided hypersurface separating two smooth initial data sets.
Suppose that the $2$-tensors $k_-$ and $k_+$ have traces along $\Sigma$ given by
\[
k_-|_{\Sigma} = \tfrac{a}{n-1} \, g_-, \qquad k_+|_{\Sigma} = \tfrac{b}{n-1} \, g_+,
\]
where $a,b \in \bR$ are constants, and the induced metrics agree on $\Sigma$.
Let $H_-$ and $H_+$ denote the mean curvatures computed on the two sides of the hypersurface $\Sigma$ with respect to a common normal.

Given a smooth positive function $t : \Sigma \to (0,\infty)$, we define the hyperbolic angle $\vartheta = - \log t$.
Then, the crease condition~\cite{KazarasKhuriLin2025}*{(1.2)} is equivalent to
\[
[t (H_- + a) - (H_+ + b) ] \cdot [ (H_-  - a) - t (H_+ - b)] \geq t^{-1} |\nabla^{\Sigma} t|^2
\]
where both factors are non-negative.
\end{lemma}
\begin{proof}
In the frame considered in~\eqref{eqn:bartnik-data}, we compute that
\[
\textup{tr}_{\Sigma} k_- = a, \qquad \textup{tr}_{\Sigma} k_+ = b, \qquad \beta_- = \beta_+ = 0.
\]
Thus, we can identify the mean curvature vectors with the vectors $\vec{H}_- = (H_-, a)$ and $\vec{H}_+ = (H_+, b)$ in the Minkowski space $\mathbb{R}^{1,1}$, equipped with the quadratic form $x_1^2-x_2^2$.
In the next step, we use the definition $t := e^{- \vartheta}$ to compute that
\begin{align*}
    X_1 + X_2 &= (\cosh \vartheta - \sinh \vartheta) H_- + (\cosh \vartheta - \sinh \vartheta)a - H_+ - b, \\
    X_1 - X_2 &= (\cosh \vartheta + \sinh \vartheta) H_- - (\sinh \vartheta + \cosh \vartheta) a - H_+ + b. 
\end{align*}
We therefore obtain
\begin{align*}
X_1^2 - X_2^2 &= [t (H_- + a) - (H_+ + b) ] \cdot [ (H_-  - a) - t (H_+ - b)]. 
\end{align*}
On the other hand, the difference of the normal connection forms appearing in~\cite{KazarasKhuriLin2025}*{\S 1} is
\[
\beta^{\Delta} = \beta_+ - \beta_- - d \vartheta = - d \vartheta = d ( \log t) = t^{-1} dt.
\]
Hence, $|\beta^{\Delta}|^2 = t^{-2} |\nabla^{\Sigma} t|^2$.
The property~\eqref{eqn:crease-condition} is equivalent to the non-negativity of both null components, and $X_1^2-X_2^2 \geq |\beta^{\Delta}|^2$.
Substituting the preceding identities proves our assertion.
\end{proof}

We note that Lemma~\ref{lemma:KKL} has several elementary consequences, which recover and extend the families sufficient conditions considered in the lock principle of~\cite{FrenckHankeHirsch2026}*{Lemma 2}.
In particular, this implies that Condition $(iii)$ is flexible and can be satisfied by a range of possible configurations of the mean curvatures of the hypersurface $\Sigma_i$, computed from the two sides.
In the following Corollary, we record a number of examples of data $(H_-, H_+, a, b, t)$ satisfying this condition.

\begin{lemma}\label{lemma:satisfy-the-condition}
Let $\Sigma$ be a smooth two-sided hypersurface separating two smooth initial data sets, and let $H_-$ and $H_+$ denote the mean curvatures computed on the two sides of the hypersurface $\Sigma$ with respect to a common normal.
Moreover, we fix two real numbers $a,b \in \mathbb{R}$.

Then, there exists a smooth positive function $t: \Sigma \to (0,\infty)$ satisfying the condition
\begin{equation}\label{eqn:t-condition}
[t (H_- + a) - (H_+ + b) ] \cdot [ (H_-  - a) - t (H_+ - b)] \geq t^{-1} |\nabla^{\Sigma} t|^2,
\end{equation}
with both factors on the left-hand side non-negative, whenever any of the following conditions holds.
\begin{enumerate}[(i)]
    \item We have $H_- - H_+ \geq |a-b|$.
    \item There exist constants $\underline{H}_-$ and $\overline{H}_+$ satisfying $\underline{H}_-^2 - a^2 \geq \overline{H}^2_+ - b^2$, as well as the properties
    \[
    H_- \geq \underline{H}_- > - a, \qquad H_+ \leq \overline{H}_+, \qquad \overline{H}_+ > -b.
    \]
    \item We have $H_- > |a|$ and $H_+ > |b|$, as well as
    \[
    \left| \nabla^{\Sigma} \log \left( \frac{H_+ + b}{H_- + a} + \frac{H_- - a}{H_+ - b} \right) \right| \leq \frac{H_-^2 - H_+^2 + b^2 - a^2}{\sqrt{2(H_-^2 + H_+^2 - a^2-b^2)}}.
    \]
    \item We have $H_- > |a|$ and $H_+ > |b|$, as well as
    \[
    \sqrt{H_-^2-a^2} - \sqrt{H_+^2-b^2} \geq \left| \frac{a \nabla^{\Sigma} H_-}{H_-^2 - a^2} - \frac{b \nabla^{\Sigma} H_+}{H_+^2-b^2} \right|.
    \]
    \item We have $H_- + a > 0$ as well as
    \[
    \sqrt{H_-^2 - a^2+b^2} > \max \{ 0, -b\}, \qquad \sqrt{H_-^2-a^2+b^2} \geq H_+ + \Bigl| \nabla^{\Sigma} \log \Bigl( \tfrac{b + \sqrt{H_-^2-a^2+b^2}}{H_-+a} \Bigr) \Bigr|.
    \]
\end{enumerate}
In the setting of Theorem~\ref{thm:better-theorem}, where $a,b \geq 0$, the assumptions $H_- > |a|$ and $H_+ > |b|$ in $(iii)$ and $(iv)$ reduce to $H_- > a$ and $H_+ > b$, respectively.
\end{lemma}
\begin{proof}
We divide the desired inequality~\eqref{eqn:t-condition} to see that it is equivalent to
\begin{equation}\label{eqn:equivalent-inequality}
    [ t (H_- + a) - (H_+ + b)] \cdot [ t^{-1}(H_- - a) - (H_+-b)] \geq |\nabla^{\Sigma} \log t|^2
\end{equation}
with both factors on the left-hand side non-negative.

In case $(i)$, the assertion follows immediately by taking $t \equiv 1$; in case $(ii)$, we choose the constant $t \equiv ( \underline{H}_- + a)^{-1} (\overline{H}_+ + b)$.
By the assumption, we obtain
\[
t(H_-+a) - (H_++b) \geq0, \qquad t^{-1} (H_- - a) - (H_+-b) \geq \frac{\underline{H}_-^2-a^2-\overline{H}^2_++b^2}{\overline{H}_++b} \geq 0.
\]
Since $\nabla^{\Sigma}t=0$, the asserted inequality~\eqref{eqn:equivalent-inequality} follows.

To treat the cases $(iii)$ and $(iv)$, we define the functions
\[
A_+ := (H_- +a)^{-1} (H_++b), \qquad A_- := (H_+-b)^{-1} (H_ - a).
\]
In both cases, our hypotheses imply that $0 < A_+ \leq A_-$.
Moreover, the non-negativity of both factors in~\eqref{eqn:equivalent-inequality} is equivalent to $A_+ \leq t \leq A_-$, while the inequality~\eqref{eqn:equivalent-inequality} itself reduces to
\begin{equation}\label{eqn:a-plus-a-minus}
(H_-+a)(H_+-b) (t-A_+)(A_--t) t^{-1} \geq |\nabla^{\Sigma} \log t|^2.
\end{equation}
In case $(iii)$, we form the arithmetic mean $t := \frac{1}{2} (A_-+A_+)$.
The left-hand side of~\eqref{eqn:a-plus-a-minus} equals
\[
\tfrac{1}{2} (H_-^2+H_+^2-a^2-b^2)^{-1} (H_-^2-H_+^2+b^2-a^2)^2.
\]
Hence, the assertion follows from the hypothesis.
In case $(iv)$, we take the geometric mean $t := \sqrt{A_-A_+}$ to compute that the left-hand side of~\eqref{eqn:equivalent-inequality} equals
\[
\bigl( \textstyle{\sqrt{H_-^2-a^2}} - \textstyle{\sqrt{H_+^2-b^2}} \bigr)^2.
\]
On the other hand, we also compute that
\[
\nabla^{\Sigma} \log t = \frac{a \nabla^{\Sigma} H_-}{H_-^2-a^2} - \frac{b \nabla^{\Sigma} H_+}{H_+^2-b^2}.
\]
Thus, the property~\eqref{eqn:equivalent-inequality} again follows from our hypotheses.

Finally, in case $(v)$, we choose
\[
t := (H_-+a)^{-1} \bigl[ b + \textstyle{\sqrt{H_-^2-a^2+b^2}} \, \bigr].
\]
A direct computation shows that
\[
t(H_-+a)-(H_++b) = t^{-1} (H_--a) - (H_+-b) = \textstyle{\sqrt{H_-^2-a^2+b^2}} - H_+.
\]
Consequently, both the non-negativity of the two factors and the inequality~\eqref{eqn:equivalent-inequality} follow directly from the hypothesis.
This completes the proof of our assertion.
\end{proof}

Next, we record a modification of the results in~\cite{KazarasKhuriLin2025}*{\S 3 - 5}, which allows us to apply their Theorem 1.4 for creased initial data in the presence of multiple corners.
See also~\cite{FrenckHankeHirsch2026}.
\begin{proposition}\label{prop:pmt-for-creased}
Suppose that Assumption~\ref{assumption} holds.
Let $k_i$ be Lipschitz symmetric two-tensors on the compact pieces $M_i$, with smooth boundary traces on the hypersurfaces $\Sigma_i$, and with $k_N =0$.
Suppose that the dominant energy condition holds almost everywhere on each piece and that the property~\eqref{lemma:KKL} is satisfied along each $\Sigma_i$, for a smooth hyperbolic angle $\vartheta_i$.
Then, $m_{\textup{ADM}}(M,g) \geq 0$ holds on every asymptotically flat end.
Moreover, if the mass vanishes on any end, then $\mu_i=0$ and $J_i=0$ hold a.e.~ on every piece.
\end{proposition}
\begin{proof}
The argument of Kazaras-Khuri-Lin~\cite{KazarasKhuriLin2025}*{Theorem 1.4} extends to our setting with minor modifications.
Their boundary computation~\cite{KazarasKhuriLin2025}*{Proposition 3.3} applies separately to each crease, while their weighted Poincar\'e inequality~\cite{KazarasKhuriLin2025}*{Corollary 4.3} follows by successive application of Lemma 4.2 therein, and the properties of the transmission boundary value problem.

These arguments continue to apply if $k \in W^{1,\infty}_{\textup{loc}}$ is only assumed to be Lipschitz.
Indeed, ~\cite{KazarasKhuriLin2025}*{Lemma 4.2} only requires the boundedness of $k$, while the integrated mass inequality of Proposition 4.4 therein follows using the fact that $k$ can be approximated uniformly and in $W^{1,p}$, for finite $p$.
Here, we also used the fact that weakly parallel spinors are locally Lipschitz, by using the parallel spinor equation and iteratively applying the Sobolev embedding theorem; see, for example, ~\cite{kapovitch-lott}*{(2.17), p. 744}.
Furthermore, in~\cite{KazarasKhuriLin2025}, the Dirac-Witten operator in the boundary regularity result of Proposition 5.1, and in the existence argument of Theorem 4.5, depends on the $2$-tensor $k$ only through a zeroth-order term.
Hence, these results can again be applied in our setting.

Applying the above discussion, we deduce that, for every constant spinor $\psi_{\infty}$ on a chosen end, there exists a Dirac-Witten harmonic spinor solving the transmission boundary value problem of~\cite{KazarasKhuriLin2025}*{\S 4} and satisfying the mass inequality
\[
\tfrac{1}{2} (n-1) \omega_{n-1} \, m_{\textup{ADM}}(M,g) \, |\psi_{\infty}|^2 \geq \sum_{i=0}^N \int_{M_i} \bigl[ |\bar{\nabla} \psi_i|^2 + \tfrac{1}{2} \la \psi_i, (\mu_i + \la J_i, \tau \rg) \psi_i \rg \bigr] \, dV.
\]
Because the dominant energy condition holds a.e.~ on each piece, we deduce that $m_{\textup{ADM}}(M,g) \geq 0$.

Finally, if $m_{\textup{ADM}}(M,g) =0$, then all spinors considered above are parallel with respect to $\bar{\nabla}$.
By the uniqueness of the parallel transport, and the invertibility of the transmission maps in~\cite{KazarasKhuriLin2025}*{\S 4}, a basis of asymptotic spinors can be transported to a basis at every point of each piece.
The mass inequality used above therefore implies that $\mu_i + \la J_i , \tau \rg = 0$ a.e.
This proves our assertion.
\end{proof}

\section{Proof of the main theorem}

We now combine the above ingredients to prove Theorem~\ref{thm:better-theorem}.

\begin{proof}[Proof of Theorem~\ref{thm:better-theorem}]
We construct an auxiliary initial data set by taking $k_i = \frac{c_i}{n-1} g_i$ on each $M_i$.
On $M_0$, we set $c_0 = a_1$; we also take $c_N = b_N=0$ on $M_N$.
First, condition $(i)$ ensures that the property~\eqref{eqn:dec-equivalent} holds on those pieces, for $c_0$ and $c_N$.

In the next step, Condition $(ii)$ allows us to apply Lemma~\ref{lemma:band-width-estimate}, for each $1 \leq i \leq N-1$.
Thus, we obtain a Lipschitz function $c_i$ on $M_i$ satisfying~\eqref{eqn:dec-equivalent}, such that $c_i = b_i$ in a neighborhood of $\Sigma_i$ and $c_i = a_{i+1}$ in a neighborhood of $\Sigma_{i+1}$.
If necessary, we may interchange the two boundary components in Lemma~\ref{lemma:band-width-estimate}; in particular, if $b_i = a_{i+1}$, we can take $c_i$ to be constant.

By our construction, and the computation~\eqref{eqn:dec-equivalent}, it follows that the auxiliary data, with $k_i = \frac{c_i}{n-1} g_i$, satisfy the dominant energy condition on every piece.
At $\Sigma_i$, the traces of the two tensors $k_{i-1}$ and $k_i$ are $\frac{a_i}{n-1} g_i$ and $\frac{b_i}{n-1} g_i$, respectively.
Moreover, Lemma~\ref{lemma:KKL} shows that condition $(iii)$ is precisely the crease condition for the angle $\vartheta_i = - \log t_i$.
Thus, we can apply the modified positive mass theorem for creased initial data of Proposition~\ref{prop:pmt-for-creased} to obtain $m_{\textup{ADM}}(M,g) \geq 0$ on each asymptotically end.

Finally, we observe that if the mass of an end vanishes, then the last assertion of Proposition~\ref{prop:pmt-for-creased} gives $d c_i = - J_i = 0$ on each piece.
Consequently, the coefficients must be constant on the connected pieces, implying that $b_i = a_{i+1}$ for all $1 \leq i < N$.
This completes the proof.
\end{proof}

In view of Lemma~\ref{lemma:satisfy-the-condition}, we obtain several concrete configurations that force non-negative ADM mass, despite the scalar curvature failing to be non-negative in the distributional sense.
For example, the conditions~\eqref{eqn:ri-condition} and~\eqref{eqn:integral-Phi} produce effective bounds in the presence of shields where $R_{g_i} \geq \kappa_i n(n-1)$, where the constant $\kappa_i$ is allowed to be negative.
Combined with cases $(i)$ and $(ii)$ of Lemma~\ref{lemma:satisfy-the-condition}, these imply positive mass theorems under $L^{\infty}$ bounds on the mean curvatures of the corners.

\bibliography{ref}

\end{document}